\documentclass{article}
\usepackage{amsmath,amsthm,amsfonts,amssymb,amscd,enumerate}

\theoremstyle{plain}
\newtheorem{theorem}{Theorem}[section]
\newtheorem{proposition}[theorem]{Proposition}

\newtheorem{lemma}[theorem]{Lemma}

\newtheorem*{Theorem}{Theorem}

\theoremstyle{definition}

        {\end{enumerate}}

\newenvironment{enumeratea'}{\begin{enumerate}[\upshape
(a)$'$]}{\end{enumerate}}

\newcommand{\comment}[1]{}

\numberwithin{equation}{section}
\begin{document}
\title{A Note on Parabolic Subgroups of a Coxeter Group}
\author{{Dongwen Qi}\\
{\small Department of Mathematics, Ohio State University}\\
{\small  231 West 18th Avenue, Columbus, OH 43210, U.S.A.}}
\maketitle
\begin{center}
\begin{minipage}{110mm}
\vskip 1cm
\begin{center}{\bf Abstract}\end{center}
{The aim of this note is to prove that the parabolic closure of
any subset of a Coxeter group is a parabolic subgroup. To obtain
that, several technical lemmas on the root system of a parabolic
subgroup are established. \vskip 0.3cm

 \vskip
0.3cm

}
\end{minipage}
\end{center}
\vskip 1cm
\baselineskip 16pt
\section{Introduction}
A Coxeter group $(W,S)$ is a group with a presentation of the
form,
\begin{equation}
\label{coxdef}
 W = {\rm gp}\langle\{s | s\in S\} | (st)^{m_{s t}}=1 \mbox{
for  all }    s, t\in S\rangle,
\end{equation}
where $m_{s t}= m_{t s}$ is a positive integer or $\infty$, and
$m_{s t}=1$ if and only if $s=t$ (A ``relation'' $(s t)^{\infty}=1$
is interpreted as vacuous). The cardinality $|S|$ of $S$ is called
the rank of $W$. The length $l(w)$ of an element $w\in W$ is the
smallest number $m$ where $w$ can be expressed as a product of $m$
elements (counting repetitions) in  $S$. We are mainly interested in
Coxeter groups of finite rank and assume $|S|$  is finite in this
note, although some statements are still valid for infinite rank
situation.

Given a Coxeter group $W$ defined as above,  for a subset $I$ of
$S$, let $W_I$ be the subgroup generated by $s\in I$ and call it a
special subgroup of $W$. At the extremes, $W_{\emptyset}=\{1\}$
and $W_S=W$. For any $w\in W$, $wW_I w^{-1}$ is called a parabolic
subgroup of $W$. The parabolic closure Pc($A$) of a subset $A$ of
$W$ is defined to be the intersection of all parabolic subgroups
containing $A$. It is believed that a parabolic closure is
a parabolic subgroup, for example, by studying the parabolic
closure of some particular element or a subgroup of a Coxeter
group, D. Krammer \cite{krammer} obtained some very interesting
results of irreducible,  infinite Coxeter groups. However, I have
not seen a proof in the literature that a parabolic closure which,
by definition, is the intersection of a collection of parabolic
subgroups, must be a parabolic subgroup. Perhaps the  result
closest to this aim is

\begin{theorem}
\label{twointer}
 The intersection of two parabolic subgroups of a
Coxeter group is a parabolic subgroup.
\end{theorem}

This result appears in geometric form in \cite{tits} and a proof
using algebraic argument is given in \cite{solomon}. However the
above proof does not establish the conclusion that a parabolic
closure
 is a parabolic subgroup.

In this note, first I give a short proof of this theorem using
standard facts of canonical representations of Coxeter groups (see
\cite{bourbaki} \cite{hump}). The proof has a simple and clear
geometric meaning. Following this, I use the general notion of
root systems developed by V. Deodhar \cite{deodroot} to establish
some technical lemmas on the root system of a parabolic subgroup
and use them to prove that the  parabolic closure  of any subset
of $W$ is a parabolic subgroup and give an alternate description
of it. All these accounts are in Section 3.

\section{Preliminaries}
In this section we collect a few basic facts about the canonical representations of
Coxeter groups. The materials  are taken from
Chapter 5 of  \cite{hump}. Let $V$
be a vector space over {\bf R}, having basis $\{\alpha_{s}| s\in
S\}$ in one-to-one correspondence with $S$. Define a symmetric
bilinear form ( , ) on $V$ by setting
\begin{equation}
\label{linearfrom} \mbox{(}\alpha_s, \alpha_t\mbox{)} =
-\cos\frac{\pi}{m_{s t}}.
\end{equation}
The value on the right-hand side is interpreted to be $ -1$ when $m_{s t}=\infty$. Now
for each $s\in S$,  define a linear transformation $\sigma_s:
V\rightarrow V$ by $\sigma_s
\lambda=\lambda-2\mbox{(}\alpha_s,\lambda\mbox{)}\alpha_s$. Then
$\sigma_s$ is an affine reflection, which has order 2 and fixes the
hyperplane $H_s=\{\delta\in V|\mbox{(}\delta,\alpha_s\mbox{)}=0\}$
pointwise,  and $\sigma_s\alpha_s=-\alpha_s$. We have

\begin{theorem}
\label{repre1}
There is a unique homomorphism $\sigma : W \rightarrow \mbox{GL(}V\mbox{)}$
sending $s$ to $\sigma_s$. This homomorphism is a faithful representation of
$W$ and the group $\sigma\mbox{(}W\mbox{)}$ preserves the bilinear form
defined as above. Moreover, for each pair $s, t\in S$, the order of $st$ in $W$ is
precisely $m_{s t}$.
\end{theorem}

Now we introduce the {\bf root system $\Phi$} of $W$, which is defined to be the
collection of all vectors $w\mbox{(}\alpha_s\mbox{)}$, where $w\in W$ and $s\in S$.
An important fact about the root system is that any root $\alpha\in \Phi$ can be
expressed as
\[\alpha=\sum\limits_{s\in S} c_s \alpha_s, \]
where all the coefficients  $c_{s}\geq 0$ (we call $\alpha$
positive and write $\alpha>0$), or all the
 coefficients $c_{s}\leq 0$ (call $\alpha$ negative and write $\alpha <0$).
 Write $\Phi^{+}$ and $\Phi^{-}$ for the respective sets of positive and negative roots.
 Then $\Phi^{+}\bigcap\Phi^{-}=\emptyset$ , $\Phi^{+}\bigcup\Phi^{-}=\Phi$
 and $\Phi^{-}=-\Phi^{+}$.  The map from $\Phi$ to $R=\{wtw^{-1}| w\in W, t\in S\}$
 (the set of reflections in $W$) given
 by $\alpha=w(\alpha_s)\mapsto  wsw^{-1}$ is well-defined
and restricts to a bijection from $\Phi^{+}$ ($\Phi^{-}$) to $R$, and $\sigma(wsw^{-1})=t_\alpha$,
where $t_{\alpha}\lambda=\lambda-2(\alpha,\lambda)\alpha$. The following
fact is important to us.

\begin{proposition}
\label{length}
Let $w\in W$, $\alpha\in\Phi^{+}$. Then $l(wt_{\alpha})>l(w)$ if and only if $w(\alpha)>0$.
\end{proposition}

With the representation $\sigma : W\rightarrow {\rm GL} (V)$ in
mind, we define a dual representation $\sigma^{*} : W\rightarrow
{\rm GL} (V^{*})$ as follows (and we abuse the notations by
identifying $w$ with $\sigma(w)$ or $\sigma^{*}(w)$),
\[\langle w(f), \lambda\rangle=\langle f, w^{-1}(\lambda)\rangle \mbox{  for  } w\in W, f\in V^{*}, \lambda\in V,\]
where $V^{*}$ is the dual space of $V$ and the natural pairing of
$V^{*}$ with $V$ is denoted by $\langle f, \lambda\rangle$. This
dual representation induces an action of $W$ on the Tits cone
defined as follows. For $I\subset S$, write
\[C_I=\{f\in V^{*}|\langle f,\alpha_s\rangle>0 \mbox{ for } s\in S-I \mbox{ and } \langle f, \alpha_s\rangle=0
\mbox{  for } s\in I\}.\] Notice that $C_S=\{0\}$ and write
$C=C_{\emptyset}$, $\overline{C}=\bigcup\limits_{I \subset S}
C_I$. Define $U$ to be the union of all $w(\overline{C})$, $w\in
W$. $U$ is a cone in $V^{*}$, called the Tits cone of $W$.

\begin{theorem}
\label{titscone} {\rm (a)} Let $w\in W$ and $I, J\subset S$. If
$w(C_I)\bigcap C_J\neq\emptyset$, then $I=J$ and $w\in W_I$, so
$w(C_I)=C_I$. In particular, $W_I$ is the precise stabilizer in
$W$ of each point of $C_I$, and $w(C_I)$, where $w\in W$,
$I\subset S$, form a partition of the Tits cone $U$.

{\rm (b)} $\overline{C}$ is a fundamental domain  for the action of
$W$ on $U$: the $W$-orbit of each point of $U$ meets $\overline{C}$
in exactly one point.

{\rm (c)} The Tits cone $U$ is convex, and every closed line
segment in $U$ meets just finitely many of the sets of the family
$\{w(C_I)|I\subset S\}$.
\end{theorem}

Both $\sigma$ and $\sigma^{*}$ are called canonical representations.

\section{Root system of a parabolic subgroup and the parabolic
closure of a set} Let me show that Theorem \ref{titscone} implies
Theorem \ref{twointer} immediately.

 \emph{Proof of Theorem \ref{twointer}}. Given two
parabolic subgroups $G_1$ and $G_2$ of $W$. Pick $x_i\in U$, $i=1,
2$, such that $G_i$ is the stabilizer of $x_i$. Then $G_1\bigcap
G_2$ fixes the line segment $\overline{x_1 x_2}$. By (c) of
Theorem \ref{titscone}, there exist $y_1\not=y_2$ on
$\overline{x_1 x_2}$ such that they belong to the same $w(C_I)$.
So $y_1$ and $y_2$ have the same stabilizer $P=wW_I w^{-1}$. Now
$P$ fixes the line segment $\overline{x_1 x_2}$ and hence
$P\subset G_i$, $P\subset G_1\bigcap G_2$. Since $G_1\bigcap G_2$
fixes $\overline{x_1 x_2}$, the reversed inclusion is obvious.
This completes the proof.

Now we describe a lemma on the root system $\Phi_I$ of a special
subgroup $W_I$, where $\Phi_I=\{w(\alpha_s)|w\in W_I, s\in I\}$.

\begin{lemma}
\label{subroot} $\Phi_I=\Phi\bigcap {\rm span}\{\alpha_s|s\in I\}$.
Here {\rm span} means {\bf R}-span.
\end{lemma}

It is obvious that $\Phi_I\subset \Phi\bigcap {\rm
span}\{\alpha_s|s\in I\}$. When $W$ is finite,  arguments similar
to that given on page 11 of \cite{hump} yields the reversed
inclusion. In the case that $W$ is of finite rank, the nontrivial
fact that $\sigma^{*}(W)$ is a discrete subgroup of GL($V^{*}$)
implies that $\Phi$ is a discrete set of $V$, which makes similar
arguments work. However, Lemma \ref{subroot} holds even when
$|S|=\infty$, as the following proof demonstrates. In fact, it
follows from the basic properties of Coxeter groups.

 \emph{Proof of Lemma \ref{subroot}}. We want to prove that
$\Phi\bigcap {\rm span}\{\alpha_s|s\in I\}\subset\Phi_I$. Pick an
arbitray $\phi\in \Phi\bigcap {\rm span}\{\alpha_s|s\in I\}$,
$\phi>0$. Write $\phi=c_1\alpha_{s_1}+\cdots +c_n\alpha_{s_n}$,
where $c_i>0$, $s_i\in I$, $i=1,\cdots, n$, $s_i\not=s_j$ when
$i\not=j$. We assume $n\geq 2$, otherwise $\phi=\alpha_{s_1}\in
\Phi_I$. Now use induction on the length $l(t_\phi)$ of $t_\phi$ .
Recall from Section 2 that
$t_\phi(\lambda)=\lambda-2(\phi,\lambda)\phi$.

Notice that $1=(\phi, \phi)=\sum\limits_{j=1}^n c_j(\phi,
\alpha_{s_j})$, we know $(\phi, \alpha_{s_i})>0$ for some $i$. A
simple calculation shows that $s_i t_\phi s_i=t_{s_i(\phi)}$ and we
want to show $l(s_i t_\phi s_i)<l(t_\phi)$. First it follows from

\begin{equation}
\label{rootsign} t_\phi(\alpha_{s_i})=\alpha_{s_i}-2(\phi,
\alpha_{s_i})\phi<0
\end{equation}

\noindent (we assume $n\geq 2$) that $l(t_\phi s_i)=l(t_\phi)-1$ by
Proposition \ref{length} and hence $l(s_i t_\phi)=l(t_\phi)-1$. If
$s_i t_{\phi}(\alpha_{s_i})>0$, then (\ref{rootsign}) implies that
$t_\phi(\alpha_{s_i})=-\alpha_{s_i}$, i.e.,
$\alpha_{s_i}-2(\phi,\alpha_{s_i})\phi=-\alpha_{s_i}$, hence
$\phi=\alpha_{s_i}$, a contradiction to the assumption $n\geq 2$.
Therefore $s_it_\phi(\alpha_{s_i})<0$ and
$l(t_{s_i(\phi)})=l(s_it_\phi s_i)=l(s_it_\phi)-1=l(t_\phi)-2$ and
induction hypothesis now applies and the proof is completed.

Now we start to discuss parabolic closures. With the notations of
Section 2, write $\Delta_K=\{\alpha_s|s\in K\}$ for $K\subset S$.
\begin{lemma}
\label{conjug} If $W_I=wW_J w^{-1}$ for some $w\in W$, $I,
J\subset S$, then $|I|= |J|$, and $w_0(\Delta_J)=\Delta_I$ for
some $w_0\in wW_J$, so $I=w_0Jw_0^{-1}$.
\end{lemma}

This lemma is stated and proved in Section 3.4 of \cite{davis1}.
The proof given there is mainly combinatorial (without using root
system), although some topological considerations (of connected
components separated by some ``walls'' of the corresponding Cayley
graph) are used. Here we give another proof.

\emph{Proof of Lemma \ref{conjug}}. We employ a few basic facts of
Coxeter groups. First, if $xtx^{-1}\in W_K$, where $x\in W$, $t\in
S$ and $K\subset S$, then $xtx^{-1}=w_1sw_1^{-1}$ for some $w_1\in
W_K$ and $s\in K$; that is, if a reflection of a Coxeter group $W$
lies in a special subgroup $W_K$, it is indeed a reflection in
$W_K$ (considering $W_K$ as a Coxeter group by itself). Second,
$wW_J=w_0W_J$, where $w_0$ satisfies that $l(w_0t)=l(w_0)+1$ for
any $t\in J$, i.e., $w_0$ is the shortest element in $wW_J$.

Now using the above $w_0$, we have $W_I=w_0W_Jw_0^{-1}$. It follows
from the correspondence of root system and reflections of Coxeter
group $W$ that $\Phi_I=w_0(\Phi_J)$. Comparing the maximal numbers
of linearly independent positive roots in these sets (By the choice
of $w_0$, $w_0(\alpha_t)>0$, for $t\in J$.), we have $|I|=|J|$.
 The fact
$\Phi_I=w_0(\Phi_J)$ implies each element of $\Phi_I$ is a positive
or negative linear combination of $w_0(\Delta_J)$, so
$\Delta_I=w_0(\Delta_J)$ and the conclusion of lemma follows.

\begin{lemma}
\label{rank} If $W_I\not\subseteq wW_J w^{-1}$, $wW_J
w^{-1}\not\subseteq W_I$, then $W_I\bigcap wW_J w^{-1}=xW_K
x^{-1}$ with $|K|<\min\{|I|,|J|\}$.
\end{lemma}

\emph{Proof of Lemma \ref{rank}}. The statement that $W_I\bigcap
wW_J w^{-1}=xW_K x^{-1}$ for some $x\in W$ and $K\subset S$ is
guaranteed by Theorem \ref{twointer}. Since $xW_K x^{-1}=x_0W_K
x_0^{-1}\subset W_I$, where $x_0$ is the shortest element in $xW_K$,
any root corresponding to a reflection in $x_0W_K x_0^{-1}$ lies in
$\Phi_I$, that is, $x_0(\Phi_K)\subset \Phi_I$. Comparing the
maximal numbers of linearly independent positive roots in these
sets, we have $|K|\leq |I|$.

Notice that $x_0(\Phi)=\Phi$, it follows from Lemma \ref{subroot}
that \[x_0(\Phi_K)=x_0(\Phi\bigcap {\rm span}\Delta_K)=\Phi\bigcap
{\rm span}\{x_0(\Delta_K)\}.\]
 If $|K|=|I|$, noticing that
$x_0(\Delta_K)\subset\Phi_I$ and ${\rm span}\{x_0(\Delta_K)\}\subset
{\rm span}\Delta_I$, we would have
\[ \Phi_I=\Phi\bigcap {\rm
span}\Delta_I=\Phi\bigcap{\rm
span}\{x_0(\Delta_K)\}=x_0(\Phi_K),\] and hence
$\Delta_I=x_0(\Delta_K)$, $W_I=x_0W_K x_0^{-1}=xW_K x^{-1}$,
contradicting the assumption of the lemma. Hence $|K|<|I|$.
Similarly, $|K|<|J|$.

Now another description of parabolic closure is

\begin{Theorem}The parabolic closure
{\rm Pc}($A$) of a subset $A$ of $W$ is the parabolic subgroup
$wW_J w^{-1}$ containing $A$, with $|J|$ being the smallest.
\end{Theorem}

The proof is obvious. The statement that the above mentioned
parabolic subgroup is contained in any parabolic subgroup containing
$A$ follows from Lemma \ref{rank} and the fact (whose proof is
essentially contained in the proof of Lemma \ref{rank}): if $xW_K
x^{-1}\subset W_I$, then $|K|\leq |I|$ and if $xW_K x^{-1}$ is a
proper subgroup of $W_I$, then $|K|<|I|$.

\end{document}